\documentclass[a4paper,12pt]{article}
\usepackage{graphicx}
\usepackage{a4wide}
\usepackage[english]{babel}
\usepackage{amsfonts}
\usepackage{amsmath}
\usepackage{amssymb}
\usepackage{dsfont}
\newtheorem{lemma}{Lemma}

\newtheorem{theorem}{Theorem}

\allowdisplaybreaks

\newcommand {\R} {\mathbb{R}}

\newcommand {\ve} {\varepsilon}

\numberwithin{equation}{section} \numberwithin{theorem}{section}
\numberwithin{lemma}{section} %\numberwithin{example}{section}
\numberwithin{corollary}{section}

\def\varep{\varepsilon}

\def\varep{\varepsilon}

\begin{document}
\title{\bf Lacunary series and stable distributions}
\author{I.\ Berkes\footnote{ Graz University of Technology, Institute of
Statistics, Kopernikusgasse 24, 8010 Graz, Austria. \mbox{e-mail}:
\texttt{berkes@tugraz.at}. Research supported by FWF grants
P24302-N18, W1230 and OTKA grant K106815.} \, and R.\
Tichy\footnote{Graz University of Technology, Institute of
Mathematics A, Steyrergasse 30, 8010 Graz, Austria. \mbox{e-mail}:
\texttt{tichy@tugraz.at}. Research supported by FWF grants
P24302-N18, W1230 and SFB project F5510.}}

\date{}
\maketitle

\begin{center}
\emph{Dedicated to Professor Paul Deheuvels on the occasion of his 65th birthday}
\end{center}

\vskip0.5cm

\abstract{By well known results of probability theory, any sequence of random variables with bounded second moments
has a subsequence satisfying the central limit theorem and the law of the iterated logarithm in a randomized form.
In this paper we give criteria for a sequence $(X_n)$ of random variables to have a subsequence $(X_{n_k})$ whose
weighted partial sums, suitably normalized, converge weakly to a stable distribution with parameter $0<\alpha<2$.}

%\vskip1cm \noindent{\bf AMS 2000 Subject classification}. Primary
%42C15, 42A55, 42A61, 30B50, 11K38, 60G50

\section{Introduction}
It is known that sufficiently thin subsequences of general r.v.\
sequences behave like i.i.d.\  sequences. For example, Chatterji
\cite{chaclt}, \cite{chalil} and Gaposhkin \cite{gap1966}, \cite{gap1972}
proved that if a sequence $(X_n)$ of r.v.'s  satisfies
$\sup_n EX_n^2<\infty$, then one can find a subsequence $(X_{n_k})$ and
r.v.'s  $X$ and $Y\ge 0$ such that
\begin{equation}\label{CLTm}
\frac{1}{\sqrt N} \sum_{k\le N} (X_{n_k} - X)
\overset{d}{\longrightarrow} N(0,Y)
\end{equation}
and
\begin{equation}\label{LILm}
\limsup_{N\to\infty} \, \frac{1}{\sqrt{2N \log\log N}} \sum_{k\le N}
(X_{n_k} - X) = Y^{1/2} \qquad \textup{a.s.},
\end{equation}
where $N(0,Y)$ denotes the distribution of the r.v.\ $Y^{1/2} \zeta$
where $\zeta$ is an $N(0, 1)$ r.v.\ independent of $Y$. Koml\'os
\cite{ko} proved that under $\sup_n E|X_n|<\infty$ there exists  a
subsequence $(X_{n_k})$ and an integrable r.v. $X$ such that
$$\lim_{N\to\infty} \frac{1}{N} \sum_{k=1}^N X_{n_k}=X \qquad \text{a.s.}$$
and Chatterji \cite{chalp} showed that under $\sup_n
E|X_n|^p<\infty$, $0<p<2$  the conclusion of the previous theorem
can be changed to
$$
\lim_{N\to\infty} \frac{1}{N^{1/p}} \sum_{k=1}^N (X_{n_k}-X)=0
\qquad \text{a.s.}
$$
for some $X$ with $E|X|^p<\infty$.
Note the randomization in all these examples:
the  role of the mean and variance of the subsequence $(X_{n_k})$ is
played by random variables $X$, $Y$. On the basis of
these and several other examples, Chatterji \cite{chasub}
formulated the following heuristic principle:

\bigskip\noindent
{\bf Subsequence Principle}. {\sl Let $T$ be a probability limit
theorem valid for all sequences of i.i.d.~random variables
belonging to an integrability class $L$ defined by the finiteness
of a norm $\|\ \cdot \|_L$. Then if $(X_n)$ is an arbitrary
(dependent) sequence of random variables satisfying
$\sup_n\|X_n\|_L < + \infty$ then there exists a subsequence
$(X_{n_k})$ satisfying $T$ in a mixed form.}
\bigskip

In a profound paper, Aldous \cite{ald} proved the validity of this principle
for all limit theorems concerning  the almost sure or distributional behavior of a
sequence of functionals $f_k(X_1, X_2, \ldots)$ of a sequence $(X_n)$ of r.v.'s.
Most "usual" limit theorems belong to this class; for precise formulations,
discussion and examples we refer to \cite{ald}. On the other hand, the theory does not
cover functionals $f_k$ containing parameters (as in  weighted
limit theorems) or allows limit theorems to involve other type of uniformities.
Such uniformities play an important role in analysis.
For example, if from a sequence $(X_n)$ of r.v.'s with finite $p$-th moments ($p\ge 1$)
one can select a subsequence $(X_{n_k})$ such that
$$
K^{-1} \left(\sum_{i=1}^N a_i^2\right)^{1/2} \le \big\Vert \sum_{i=1}^N a_i
X_{n_i}\big\Vert_p  \le  K \left(\sum_{i=1}^N a_i^2\right)^{1/2}
$$
for some constant $0<K<\infty$, for every $N\ge 1$ and every $(a_1,\ldots,a_N)\in \R^N$, then
the subspace of $L^p$ spanned by $(X_n)$ contains a subspace isomorphic to Hilbert space. Such
embedding arguments go back to the classical paper of Kadec and Pelczynski \cite{kape} and play
an important role in Banach space theory, see e.g.\ Dacunha-Castelle and Krivine \cite{dc2},
Aldous \cite{ald2}. In the theory of orthogonal series and in Banach space theory we frequently
need subsequences $(f_{n_k})$ of a sequence $(f_n)$ such that $\sum_{k=1}^\infty a_k f_{n_k}$
converges a.e.\ or in norm, after any permutation of its terms, for a class of coefficient sequences
$(a_k)$. Here we need uniformity both over a class of coefficient sequences $(a_k)$ and over all
permutations of the terms of the series. A number of uniform limit theorems for subsequences have
been proved by ad hoc arguments. R\'ev\'esz \cite{re} showed that for any sequence $(X_n)$ of r.v.'s
satisfying $\sup_n EX_n^2<\infty$ one can find a subsequence $(X_{n_k})$ and a r.v.\ $X$ such that
$\sum_{k=1}^\infty a_k (X_{n_k}-X)$ converges a.s.\ provided $\sum_{k=1}^\infty a_k^2<\infty$.
Under $\sup_n \|X_n\|_\infty <+\infty$, Gaposhkin \cite{gap1966} showed that there exists a subsequence
$(X_{n_k})$ and r.v.'s $X$ and $Y\ge 0$ such that for any real sequence $(a_k)$ satisfying the uniform
asymptotic negligibility condition
\begin{equation}\label{uan1}
\max_{1\le k \le N} |a_k|=o(A_N), \qquad A_N=\left(\sum_{k=1}^N a_k^2\right)^{1/2}
\end{equation}
we have
\begin{equation}\label{CLTm1}
\frac{1}{A_N} \sum_{k\le N} a_k (X_{n_k} - X)
\overset{d}{\longrightarrow} N(0,Y)
\end{equation}
and for any real sequence $(a_k)$ satisfying the Kolmogorov condition
\begin{equation}\label{uan2}
\max_{1\le k \le N} |a_k|=o(A_N/(\log\log A_N)^{1/2})
\end{equation}
we have
\begin{equation}\label{CLTm2}
\frac{1}{(2A_N\log\log A_N)^{1/2}} \sum_{k\le N} a_k (X_{n_k} - X)=Y^{1/2} \qquad \text{a.s.}
\end{equation}
For a fixed coefficient sequence $(a_k)$ the above results follow from Aldous' general theorems,
but the subsequence $(X_{n_k})$ provided  by the proofs depends on $(a_k)$ and to find a
subsequence working for all $(a_k)$ simultaneously requires a uniformity which is, in general,
not easy to establish and it can fail in important situations. (See Guerre and Raynaud \cite{gura}
for a natural problem where uniformity is not valid.) In \cite{ald}, Aldous used an
equicontinuity argument to prove a permutation-invariant version of the theorem of
R\'ev\'esz above, implying that every orthonormal system $(f_n)$ contains a subsequence $(f_{n_k})$
which, using the standard terminology, is an {\it unconditional convergence system.}
This had been a long standing open problem in the theory of orthogonal series (see  Uljanov \cite{ul}, p.\ 48)
and was
first proved by Koml\'os \cite{ko2}.  In \cite{be1985} we used the method of Aldous to prove
extensions of the Kadec-Pelczynski theorem, as well as selection theorems for almost symmetric sequences.
The purpose of the present paper is to use a similar technique to prove a uniform limit theorem of
probabilistic importance, namely the analogue of Gaposhkin's uniform CLT (\ref{uan1})--(\ref{CLTm1})
in the case when the limit distribution of the normed sum is a stable law with parameter $0<\alpha<2$.
To formulate our result, we need some definitions. Using the terminology of \cite{bero}, call the
sequence $(X_n)$ of r.v.'s {\it determining} if it has a limit distribution relative to any set $A$
in the probability space with $P(A) > 0$, i.e.~for any $A\subset\Omega$ with $P(A) > 0$ there exists
a distribution function $F_A$ such that
$$\lim\limits_{n \to\infty} P(X_n < t\mid A) = F_A(t)$$
for all continuity points $t$ of $F_A$. By an extension of the
Helly-Bray theorem (see \cite{bero}), every tight sequence of r.v.'s contains
a determining subsequence. Hence in studying the asymptotic behavior of thin subsequences of general
tight sequences we can assume without loss of generality that our original sequence $(X_n)$ is
determining. By \cite{bero}, Proposition 2.1, for any continuity point $t$ of the limit distribution
function $F_\Omega$, the sequence $I\{X_n\le t\}$ converges weakly in $L^\infty$ to some r.v. $G_t$;
clearly $G_s \le G_t$ a.s.\ for any $s\le t$. (A  sequence $(\xi_n)$ of bounded r.v.'s is said
to converge to a bounded r.v.\ $\xi$ weakly in $L^\infty$ if $E(\xi_n\eta)\longrightarrow E(\xi \eta)$ for any
integrable r.v.\ $\eta$. To avoid confusion, we will call ordinary weak convergence of probability theory
distributional convergence). Using a standard procedure (see e.g.\ R\'ev\'esz \cite{re2}, Lemma 6.1.4), by choosing a dense
countable set $D$ of continuity points of $F_\Omega$, one can construct versions of $G_t$, $t\in D$ such
that, for every fixed $\omega\in \Omega$, the function $G_t (\omega), t\in D$ extends to a distribution function.
Letting  $\mu$ denote the corresponding measure,  $\mu$ is called the {\it limit random measure\/} of $(X_n)$;
it was introduced by Aldous \cite{ald}; for properties and applications see \cite{ald2}, \cite{be1985},
\cite{bepe}, \cite{bero}. Clearly, $\mu$ can be considered as a measurable map from the underlying probability
space $(\Omega,{\cal F}, P)$ to the space $\mathcal{M}$ of probability
measures on $\mathbb{R}$ equipped with the Prohorov metric $\pi$. It is easily seen that for any $A$ with
$P(A) > 0$ and any continuity point $t$ of $F_A$
we have
\begin{equation}\label{(4)} F_A(t) =
E_A(\mu (-\infty, t)),
\end{equation}
where $E_A$ denotes conditional expectation given $A$. Note that $\mu$ depends on the actual
r.v.'s $X_n$, but the distribution of $\mu$ in $(\mathcal{M}, \pi)$ depends solely on the distribution
of the sequence $(X_n)$.
The situation  concerning the unweighted CLT for lacunary sequences can now be summarized by the following
theorem.
\begin{theorem} \label{theorem1} Let $(X_n)$ be a determining sequence of r.v.'s with limit
random measure $\mu$. Then there exists a subsequence $(X_{n_k})$ satisfying, together with
all of its subsequences, the CLT (\ref{CLTm}) with suitable r.v.'s $X$  and $Y\ge 0$ if and
only if
\begin{equation}\label{cltcond}
\int_{-\infty}^\infty x^2 d\mu(x) < \infty \qquad \text{a.s.}
\end{equation}
\end{theorem}

The sufficiency part of the theorem is contained in  Aldous'general subsequence theorems in \cite{ald}; the necessity
was proved in our recent paper \cite{bt}.  Note that the condition for the CLT for lacunary subsequences of $(X_n)$
is given in terms of the limit random measure of $(X_n)$ and this condition is the exact analogue of the condition
in the i.i.d.\ case, only the common distribution of the i.i.d.\ variables is replaced by the limit random measure.
%namely the existence of finite second moments of $\mu$ for almost all $\omega$.
Note also that the existence of second moments of $(X_n)$ (or the existence of any moments)
is not necessary for the conclusion of Theorem \ref{theorem1}.

\bigskip
In this paper we investigate the analogous question in case of a nonnormal stable limit distribution,
i.e.\ the question under what conditions a sequence $(X_n)$ of r.v.'s
has a subsequence $(X_{n_k})$ whose weighted partial sums, suitably normalized, converge weakly to
an $\alpha$-stable distribution, $0<\alpha<2$.
Let, for $c>0$ and $0<\alpha<2$, $G_{\alpha, c}$ denote the distribution function with characteristic function
$\exp(-c|t|^\alpha)$ and let $S=S(\alpha, c)$ denote the class of symmetric distributions on $\mathbb R$
with characteristic function $\varphi$ satisfying
\begin{equation}\label{phi}
\varphi(t)=1-c|t|^\alpha+o(|t|^\alpha) \qquad \text {as} \ t\to 0.
\end{equation}
Our main result is

\begin{theorem}\label{theorem2}   Let $0<\alpha<2$, $c>0$  and let $(X_n)$ be a determining sequence of
r.v.'s with limit random measure $\mu$.
Assume that $\mu\in S(\alpha, c)$ with probability 1. Then there exists a subsequence $(X_{n_k})$
such that for any real sequence $(a_k)$ satisfying
\begin{equation}\label{an}
\max_{1\le k \le N} |a_k|=o(A_N), \quad A_N=\left(\sum_{k=1}^N |a_k|^\alpha \right)^{1/\alpha}
\end{equation}
we have
$$A_N^{-1} \sum_{k=1}^N a_k X_{n_k}\overset{d}{\longrightarrow} G_{\alpha, c}.$$
\end{theorem}

Condition (\ref{phi}) holds provided the corresponding (symmetric) distribution function $F$ satisfies
$$ 1-F(x)= c_1 x^{-\alpha} + \beta (x) x^{-\alpha}, \qquad x>0$$
where $c_1>0$ is a suitable constant, $\beta (x)$ is non-increasing for $x\ge x_0$ and $\lim_{x\to\infty} \beta (x) = 0$.
(See Berkes and Dehling \cite{bede}, Lemma 3.2.) Apart from the monotonicity condition, this is equivalent to the fact
that $F$ is in the domain of normal attraction of a symmetric stable distribution. (See e.g.\ Feller \cite{fe}, p.\ 581.)
It is natural to ask if the conclusion of
Theorem \ref{theorem2} remains valid (with a suitable centering factor) assuming only that $\mu \in S$
a.s.\ where $S$ denotes the domain of normal attraction of a fixed stable distribution.  From the theory
in \cite{ald} it follows that the answer is affirmative in the unweighted case $a_k=1$, but in the uniform
weighted case the question remains open. Symmetry plays no essential role in the proof of Theorem \ref{theorem2};
it is used only in Lemma \ref{lemma1}
and at the cost of minor changes in the proof, (\ref{phi}) can be replaced by a condition covering nonsymmetric distributions
as well. But since we do not know the optimal condition, we restricted our investigations to the case (\ref{phi})
where the technical details are the simplest and the idea of the proof becomes the most transparent.

Given a sequence  $(X_n^*)$ of r.v.'s and a random measure $\mu$ defined on a probability space $(\Omega, \mathcal{F}, P)$
such that $X_n^*$ are conditionally i.i.d.\ given $\mu$ with conditional distribution $\mu$, the limit random measure of
$(X_n^*)$ is easily seen to be $\mu$. The sequence $(X_n^*)$ is exchangeable, so passing to subsequences does not change
its asymptotic properties, so if $\mu \in S(\alpha, c)$ a.s., then the conclusion of Theorem \ref{theorem2} holds for the
whole sequence $(X_n^*)$ without passing to any subsequence. (This follows directly also from Lemma \ref{lemma1}.)
%This example is typical in the sense that
%provides a simplest special case for Theorem \ref{theorem2}; the theorem states
Theorem \ref{theorem2} shows that any deterministic sequence $(X_n)$ with a limit random measure $\mu$ satisfying   $\mu \in S(\alpha, c)$
a.s.\ has a subsequence $(X_{n_k})$ whose weighted partial sums behave, in a uniform sense, similarly to those of $(X_n^*)$.

\section{Proof of Theorem \ref{theorem2}}

As the first  step of the proof, we select a sequence $n_1<n_2<\ldots$ of integers such that, after a suitable discretization
of $(X_n)$, we have
\begin{equation}\label{condconv}
 P(X_{n_k} \in J |X_{n_1}, \ldots, X_{n_{k-1}})(\omega) \longrightarrow \mu (\omega, J) \quad \text{a.s.}
\end{equation}
for a large class of intervals $J$.  This step follows exactly Aldous \cite{ald}, see Proposition 11 of \cite{ald}
for details.
%Relation (\ref{condconv}) means that the conditional distributions of $(X_{n_k})$ resemble those of
%an exchangeable sequence with conditional distribution $\mu$.
Let $(Y_n)$ be a sequence of
r.v.'s on $(\Omega ,{\cal F}, P)$ such that, given {\bf X} and $\mu$, the r.v.'s $Y_1,Y_2,\ldots\ $ are conditionally
i.i.d.\ with distribution $\mu$, i.e.,
\begin{equation}\label{18}
P(Y_1\in B_1,\ldots,Y_k\in B_k\vert {\bf X},\mu) = \prod_
{i=1}^k P(Y_i\in B_i\vert {\bf X},\mu) \ \hbox{ a.s.}
\end{equation}
\begin{equation}\label{19}
P(Y_j\in B\vert {\bf X},\mu )= \mu (B)\ \hbox{ a.s.}
\end{equation}
for any $j,k$ and Borel sets $B,B_1,\ldots,B_k$ on the real line. Such
a sequence $(Y_n)$ always exists after redefining $(X_n)$ and $\mu$ on a suitable,
larger probability space; for example, one can define the triple $((X_n), \mu, (Y_n))$
on the product space $\mathbb{R}^\infty \times \mathcal{M} \times  \mathbb{R}^\infty$
as done in \cite{ald}, p.\ 72.  This redefinition will not change the distribution
of the sequence $(X_n)$ and thus by Proposition 2.1 of \cite{bero} it remains determining.
Since the random measure $\mu$ depends on the variables $X_n$ themselves and not only on
the distribution of $(X_n)$, this redefinition will change $\mu$, but not the joint
distribution of  $(X_n)$ and $\mu$ on which our results depend. Using (\ref{condconv}) and
a martingale argument, in \cite{ald}, Lemma 12 it is shown that

\begin{lemma}\label{lemma2} For every $\sigma ({\bf X})$-measurable {\rm r.v.}\
$Z$ and any $j\ge 1$ we have
$$(X_{n_k}, Z) \buildrel {d}\over \longrightarrow (Y_j,Z) \quad \text{as}\  k\to \infty.$$
\end{lemma}

\medskip
We now construct a further subsequence of $(X_{n_k})$ satisfying the conclusion of Theorem \ref{theorem2}.
By reindexing our variables, we can assume that Lemma \ref{lemma2} holds with $n_k=k$. For our construction
we need some auxiliary considerations. For a (nonrandom) measure $\mu\in S(\alpha, c)$, the
corresponding characteristic function $\varphi$ satisfies
\begin{equation}\label{sac}
\varphi(t)=1-c|t|^\alpha  + \beta (t) |t|^\alpha, \qquad t \in \mathbb{R}
\end{equation}
where $\beta$ is a bounded continuous function on $\mathbb R$ with $\beta (0)=0$. Given $\mu_1, \mu_2 \in S(\alpha, c)$
with characteristic functions $\varphi_1, \varphi_2$ and corresponding functions $\beta_1, \beta_2$ in (\ref{sac}),
%corresponding to $\varphi_1, \varphi_2$. ,
define
\begin{equation}\label{beta}
\rho (\mu_1, \mu_2)=\sup_{0\le |t|\le 1} |\beta_1(t)-\beta_2 (t)|+ \sum_{k=0}^\infty \frac{1}{2^k} \sup_{2^k \le |t|\le 2^{k+1}} |\beta_1(t)-\beta_2 (t)|.
\end{equation}
 Clearly, $\rho$ satisfies the triangle inequality
and if $\rho (\mu_1, \mu_2)=0$, then $\varphi_1(t)=\varphi_2(t)$ for all $t\in\mathbb{R}$ and thus $\mu_1=\mu_2$. Hence, $\rho$ is a metric on $S(\alpha, c)$. If $\mu, \mu_1, \mu_2, \ldots \in S(\alpha, c)$ with corresponding characteristic functions $\varphi, \varphi_1, \varphi_2, \ldots$
and  functions $\beta, \beta_1, \beta_2, \ldots$, then $\rho (\mu_n, \mu)\to 0$ implies that $\beta_n (t) \to \beta (t)$ and consequently $\varphi_n (t) \to \varphi (t)$ uniformly on compact  intervals and thus $\mu_n\overset{d}{\to} \mu$. Conversely, if $\mu_n\overset{d}{\to} \mu$,
then $\varphi_n (t) \to \varphi (t)$ uniformly on compact  intervals and thus $\beta_n (t) \to \beta (t)$  uniformly on compact intervals not containing 0. Note that $\lim_{t\to 0} \beta_n (t)=0$ for any fixed $n$ by the definition of $S(\alpha, c)$;  if this relation holds uniformly
in $n$, then $\beta_n (t) \to \beta (t)$  will hold uniformly also on all compact intervals containing 0 and upon observing that (\ref{sac})
implies $ |\beta(t)|\le |t|^{-\alpha} |\varphi(t)-1|+c \le c+2$ for $|t|\ge 1$ and thus the total contribution of the terms of the sum in
(\ref{beta}) for $k\ge M$ is $\le 4(c+2)2^{-M}$, it follows that $\rho (\mu_n, \mu)\to 0$. Thus if for a class
$H\subset S(\alpha, c)$ we have $\lim_{t\to 0} \beta (t)=0$ uniformly for all  functions $\beta$ corresponding to
measures in $H$, then in $H$ convergence of elements in Prohorov metric and in the metric $\rho$ are equivalent.

Let now $\varphi (t)=\varphi (t, \omega)$ denote the characteristic function of the random measure $\mu=\mu (\omega)$.
By the assumption $\mu\in S(\alpha, c)$ a.s.\ of Theorem \ref{theorem2}, we have
\begin{equation} \label{phiomega}
\varphi (t, \omega)=1-c|t|^\alpha +\beta(t, \omega) |t|^\alpha, \qquad t\in {\mathbb R}, \ \omega\in \Omega
\end{equation}
where $\lim_{t \to 0} \beta(t, \omega)=0$ a.s. Let $\xi_n (\omega)=\sup_{|t|\le 1/n} |\beta(t, \omega)|$, then
$\lim_{n\to\infty} \xi_n(\omega)=0$ a.s.\ and thus by Egorov's theorem (see \cite{eg}) for any $\ve>0$ there exists a measurable
set $A\subset \Omega$ with $P(A)\ge 1-\ve$  such that $\lim_{n\to\infty} \xi_n(\omega)=0$ and consequently $\lim_{t \to 0}
\beta(t, \omega)=0$ uniformly on $A$.
%, i.e.\ there
%exists a numerical (nonrandom) sequence $\ve_n\to 0$ such that $|\xi_n(\omega)|\le \ve_n$ for $\omega\in A, n=1, 2, \ldots$.
%Letting $\beta^*(t)=\ve_n$ for $1/(n+1)<t\le 1/n$, $n=1, 2, \ldots$, we have
%$$|\beta(t, \omega)| \le \beta^*(t), \quad  0<t\le 1, \, t\in A.$$
Considering $A$ as a new probability space, we will show that there exists a subsequence $(X_{n_k})$ (depending on $A$)
satisfying the conclusion of Theorem \ref{theorem2} together with all its subsequences. By a diagonal argument we can get
then a subsequence $(X_{n_k})$ satisfying the conclusion of Theorem \ref{theorem2} on the original $\Omega$.
Thus without loss of generality we can assume in the sequel that the function $\beta(t, \omega)$ in (\ref{phiomega}) satisfies
$\lim_{t\to 0} \beta(t, \omega)=0$ uniformly in  $\omega\in \Omega$ and thus by the remarks in the previous paragraph, in the
support of the random measure $\mu$ the Prohorov metric and the metric $\rho$ generate the same convergence.

\begin{lemma}\label{lemma1} Let $\mu_1, \mu_2 \in S(\alpha, c)$ satisfying (\ref{phi}), let $Z_1, \ldots, Z_n$ and $Z_1^*, \ldots, Z_n^*$
be i.i.d.\ sequences with respective distributions $\mu_1$, $\mu_2$. Let $(a_1, \ldots, a_n) \in {\mathbb R}^n$,
$A_n=\left(\sum_{k=1}^n |a_k|^\alpha\right)^{1/\alpha}$,  $\delta_n=\max_{1\le k\le n} |a_k|/A_n$. Then for $|t|\delta_n\le 1$ we have
\begin{equation}\label{chf}
\left|E\exp \left( it A_n^{-1} \sum_{k=1}^n a_k Z_k\right)-E\exp \left( it A_n^{-1} \sum_{k=1}^n a_k Z_k^*\right)\right|
\le |t|^{\alpha} \rho(\mu_1, \mu_2)
\end{equation}
where $\rho$ is defined by (\ref{beta}).
\end{lemma}

\bigskip\noindent
{\bf Proof.} Letting $\varphi_1$, $\varphi_2$ denote the characteristic function of the $Z_k$'s resp. $Z_k^*$'s and
using (\ref{sac}), (\ref{an}) and the inequality
$$ \left|\prod_{k=1}^n x_k-\prod_{k=1}^n y_k\right|\le \sum_{k=1}^n |x_k-y_k|,$$
valid for all $|x_k|\le 1, |y_k|\le 1$ we get that for $|t|\delta_n\le 1$ the left hand side of  (\ref{chf}) equals
\begin{align*}
&\left| \prod_{k=1}^n \varphi_1( ta_k/A_n)- \prod_{k=1}^n \varphi_2( ta_k/A_n)\right|\le
\sum_{k=1}^n \left|\varphi_1( ta_k/A_n)- \varphi_2( ta_k/A_n)\right|\\
&\le \sum_{k=1}^n |\beta_1 (ta_k/A_n)-\beta_2 (ta_k/A_n)| |ta_k/A_n|^\alpha \le \sup_{|x|\le |t| \delta_n }| \beta_1(x)-\beta_2 (x)|\sum_{k=1}^n |ta_k/A_n|^\alpha\\
&= |t|^\alpha \sup_{|x|\le |t|\delta_n}| \beta_1(x)-\beta_2 (x)| \le |t|^\alpha \rho(\mu_1, \mu_2).
\end{align*}

\noindent
{\bf Remark.} The proof of Lemma \ref{lemma1} shows that for any $t\in {\mathbb R}$ the left hand side of (\ref{chf}) cannot exceed
$|t|^\alpha \sup_{|x|\le |t|\delta_n}| \beta_1(x)-\beta_2 (x)|$, a fact that will be useful in the sequel.

%\bigskip\noindent
%\begin{lemma}\label{lemma2} For every $\sigma ({\bf X})$-measurable {\rm r.v.}\
%$Z$ and any $j\ge 1$ we have
%$$(X_n,Z) \buildrel {d}\over \longrightarrow (Y_j,Z).$$
%\end{lemma}

\bigskip
Given probability measures $\nu_n,\nu$ on the Borel sets of a separable
metric space $(S,d)$ we say, as usual, that $\nu_n \buildrel {d}\over
\longrightarrow \nu$ if
\begin{equation}\label{20}
\int_S f(x) d\nu_n(x) \longrightarrow \int_S f(x)
d\nu (x) \ \hbox{ as }\ n\to \infty
\end{equation}
for every bounded, real valued continuous function $f$ on $S$.
(\ref{20}) is clearly equivalent to
\begin{equation}\label{21}
Ef(Z_n) \longrightarrow Ef(Z)
\end{equation}
where $Z_n,Z$ are r.v.'s valued in $(S,d)$ (i.e.\ measurable maps
from some probability space to $(S,d)$) with distribution $\nu_n,\nu$.

\bigskip\noindent
\begin{lemma}\label{lemma3} (see \cite{rr}). Let $(S,d)$ be a separable
metric space and let $\nu,\nu_1, \nu_2, \ldots$ be probability measures
on the Borel sets of $(S,d)$ such that $\nu_n\buildrel {d}\over
\longrightarrow \nu$. Let $\cal G$ be a class of real valued functions on
$(S,d)$ such that

\medskip\noindent
(a)  $\cal G$ is locally equicontinuous, i.e.\ for for every $\varep >0$ and $ x\in S$ there is
a $\delta = \delta (\varep, x) >0$ such that $y\in S$, $d(x,
y) \le \delta$ imply $\vert f(x)-f(y)\vert \le \varep$ for
every $f\in {\cal G}$.

\smallskip\noindent
(b) There exists a continuous function $g\ge 0$ on $S$ such that
$\vert f(x)\vert \le g(x)$ for all $f\in {\cal G}$ and $ x
\in S$
\smallskip\noindent
and
\begin{equation}\label{22}
\int_S g(x) d\nu_n (x) \longrightarrow \int_S g(x) d\nu(x)\ (<\infty)
\ \hbox{ as }\ n\to \infty.
\end{equation}
Then
\begin{equation}\label{23}
\int_S f(x) d\nu_n(x) \longrightarrow \int_S f(x) d\nu(x) \ \hbox{ as }
\ n\to \infty
\end{equation}
uniformly in $f\in {\cal G}$.
\end{lemma}

Assume now that $(X_n)$ satisfies the assumptions of Theorem \ref{theorem2}, fix $t\in \mathbb{R}$
and for any $n\ge 1$, $(a_1, \ldots, a_n)\in {\mathbb R}^n$ let
\begin{equation}\label{psidef}
\psi (a_1, \ldots, a_n)=E \exp \left(itA_n^{-1} \sum_{k=1}^n a_kY_k\right),
\end{equation}
where $A_n=(\sum_{k=1}^n |a_k|^\alpha)^{1/\alpha}$ and $(Y_k)$ is the sequence of r.v.'s defined before Lemma \ref{lemma2}.
We show that for any $\varepsilon>0$ there exists a sequence $n_1<n_2<\cdots$
of integers such that
\begin{equation}\label{333}
(1-\ve) \psi (a_1,\ldots,a_k)\le E\exp \left(itA_k^{-1} \sum_{i=1}^k a_i X_{n_i}\right)\le (1+\ve) \psi (a_1,\ldots,a_k)
\end{equation}
for all $k\ge 1$ and all $(a_k)$ satisfying (\ref{an}); moreover, (\ref{333}) remains valid for every further subsequence of $(X_{n_k})$
as well.
%To begin with, let us recall that by the crucial assumption of Theorem
%\ref{theorem2}, the limit random measure $\mu$ belongs to $S(\alpha, c)$ with probability 1, i.e.\ for almost all $\omega$ the
%characteristic function of $\mu$ satisfies (\ref{phi}). The $o(|t|^\alpha)$ in (\ref{phi}) depends on $\omega$, but by a
%standard measure theoretic argument, for any $\delta>0$ there exists a set $\Omega' \subset \Omega$ with probability $\ge 1-\delta$
%such that the $o(|t|^\alpha)$ in (\ref{phi}) is uniform for $\omega \in \Omega'$. Thus using a diagonal argument we can assume, without
%loss of generality, that the characteristic function of the limit random measure $\mu$ satisfies (\ref{phi}) uniformly.
To construct $n_1$ we set
\begin{align*}
Q({\bf a},n,\ell) &= \exp\left( itA_\ell^{-1} (a_1X_n + a_2Y_2+\cdots + a_\ell
Y_\ell)\right) \cr
R({\bf a},\ell) &= \exp \left( itA_\ell^{-1} (a_1Y_1+a_2Y_2+\cdots + a_\ell Y_\ell)\right)
\end{align*}
for every $n\ge 1$, $\ell \ge 2$ and ${\bf a} = (a_1,\ldots,a_\ell) \in
R^\ell$.  We show that
\begin{equation}\label{conv}
E\left\{ {{Q({\bf a},n,\ell)}\over {\psi ({\bf a})}}
\right\} \longrightarrow  E\left\{ {{R({\bf a},\ell)}\over {\psi ({\bf a})}}
\right\} \ \hbox{ as }\ n\to \infty \quad \hbox{uniformly in }\ {\bf a},\ell .
\end{equation}
(The right side of (\ref{conv}) equals 1.)  To this end we recall that, given
{\bf X} and $\mu$, the r.v.'s $Y_1,Y_2,\ldots \ $ are conditionally
i.i.d.\ with common conditional distribution $\mu$ and thus, given
${\bf X},\mu$ and $Y_1$, the r.v.'s $Y_2,Y_3,\ldots \ $ are conditionally
i.i.d.\ with distribution $\mu$.  Thus
\begin{equation}\label{Q}
E\bigl( Q({\bf a},n,\ell)\vert {\bf X},\mu\bigr) =g^{{\bf a},\ell} (X_n,\mu)
\end{equation}
and
\begin{equation}\label{R}
E\bigl( R({\bf a},\ell)\vert {\bf X},\mu,Y_1\bigr) = g^{{\bf a},\ell}
(Y_1,\mu),
\end{equation}
where
$$g^{{\bf a},\ell} (u,\nu) = E\exp \left(itA_\ell^{-1} \left(a_1 u+\sum_{i=2}^\ell a_i\xi_i^
{(\nu)}\right) \right)  \qquad (u\in \mathbb{R}^1\ ,\ \nu \in S)$$
and $(\xi_n^{(\nu)})$ is an i.i.d.\ sequence with distribution $\nu$.
Integrating (\ref{Q}) and (\ref{R}), we get
\begin{equation}\label{26}
E\bigl( Q({\bf a}, n,\ell)\bigr) = Eg^{{\bf a},\ell}
(X_n,\mu)
\end{equation}
\begin{equation}\label{27}
E\bigl( R({\bf a},\ell)\bigr) = Eg^{{\bf a},\ell} (Y_1,\mu)
\end{equation}
and thus (\ref{conv}) is equivalent to
\begin{equation}\label{28}
E {{ g^{{\bf a},\ell} (X_n,\mu)}\over {\psi ({\bf a})}}
\longrightarrow
E {{g^{{\bf a},\ell} (Y_1,\mu)} \over {\psi ({\bf a})}} \ \hbox{ as }\
n\to \infty,\
\hbox{ uniformly in } \ {\bf a},\ell .
\end{equation}
We shall derive (\ref{28}) from Lemmas \ref{lemma2}-- \ref{lemma3}. Recall that $\rho$ is a metric
on $S = S(\alpha, c)$; the remarks at the beginning of this section show that on the support of $\mu$ the
metric $\rho$ and the Prohorov metric $\pi$ induce the same convergence and thus the same
%clearly convergence in this metric implies ordinary weak convergence
%of probability measures, i.e.\ convergence in the Prohorov metric $\pi$.
Borel $\sigma$-field; thus the limit random measure $\mu$, which is a random variable taking values in
$(S,\pi)$,
%({\it i.e.},a measurable map from the underlying probability space to $(S,{\cal B}_
%{d_0})$ where ${\cal B}_{d_0}$ denotes the Borel $\sigma$-field in $S$
%generated by $d_0$)
can be also regarded as a random variable taking values in $(S,\rho)$. Also, $\mu$ is
clearly $\sigma ({\bf X})$ measurable and thus $(X_n,\mu) \buildrel {d}\over\longrightarrow (Y_1,\mu)$
by Lemma \ref{lemma2}.  Hence, (\ref{28}) will follow from Lemma \ref{lemma3} (note the equivalence
of (\ref{20}) and (\ref{21})) if we show that the class of functions
\begin{equation}\label{29}
\left\{ {{g^{{\bf a},\ell}(t,\nu)}\over {\psi ({\bf a}) }}
\right\}
\end{equation}
defined on the product metric space $(\mathbb{R}\times S\ ,\ \lambda\times \rho)$
($\lambda$ denotes the ordinary distance on $\mathbb{R}$) satisfies conditions
(a),(b) of Lemma \ref{lemma3}. To see the validity of (a) let us note that by (\ref{18}), (\ref{19}),
$Y_n$ are conditionally i.i.d.\ with respect to $\mu$ with conditional distribution $\mu$, moreover,
we assumed without loss of generality that the characteristic function $\varphi(t, \omega)$ of
$\mu (\omega)$ satisfies (\ref{phiomega}) with $\lim_{t\to 0} \beta(t, \omega)=0$ uniformly in
$\omega$ and thus applying Lemma \ref{lemma1} with $\varphi_1(t)=\varphi(t, \omega)$ and
$\varphi_2(t)=\exp (-c|t|^\alpha)$ and using (\ref{an}) and the remark
after the proof of Lemma \ref{lemma1} it follows that there
exists an integer $n_0$ and a positive constant $c_0$ such that $\psi({\bf a})\ge c_0$ for $n\ge n_0$ and all
$(a_k)$.
Thus the validity of (a) follows from Lemma \ref{lemma1};  the validity of (b)
is immediate from $|g^{{\bf a},\ell} (u,\nu)|\le 1$. We thus
proved relation (\ref{28}) and thus also (\ref{conv}), whence it follows (note again
that the right side of (\ref{conv}) equals 1) that
\begin{equation}\label{34}
\psi({\bf a})^{-1} E\exp\left( itA_\ell^{-1} (a_1X_n + a_2Y_2+\cdots + a_\ell
Y_\ell)\right)\longrightarrow 1
\end{equation}
as $n\to \infty$, uniformly in $\ell, {\bf a}$. Hence given $\ve>0$, we can choose $n_1$ so large that
\begin{align}\label{ind1}
&|E\exp\left( itA_\ell^{-1} (a_1X_{n} + a_2Y_2+\cdots + a_\ell
Y_\ell)\right)
-E\exp(itA_\ell^{-1} (a_1 Y_1 + a_2 Y_2+\cdots + a_\ell
Y_\ell))|\nonumber\\
&\le {\varep \over 2} \psi (a_1,\ldots,a_\ell)
\end{align}
for every $\ell,{\bf a}$ and $n\ge n_1$. This completes the first induction step.

Assume now that $n_1,\ldots ,n_{k-1}$ have already been chosen.
Exactly in the same way as we proved (\ref{34}), it follows that for $\ell >k$
\begin{align*}
&\psi ({\bf a})^{-1}E \exp \left(it A_\ell^{-1} (a_1 X_{n_1} +\cdots + a_{k-1}X_{n_{k-1}}
+a_k X_n +a_{k+1}Y_{k+1}+\cdots + a_\ell Y_\ell) \right) \cr
&\longrightarrow \psi ({\bf a})^{-1}E \exp\left( itA_\ell^{-1} (a_1X_{n_1}+\cdots + a_{k-1}
X_{n_{k-1}} +a_kY_k+\cdots +a_\ell Y_\ell)\right)\ \hbox{ as }\
n\to \infty
\end{align*}
uniformly in {\bf a} and $\ell$.  Hence we can choose $n_k>n_{k-1}$ so large that
\begin{align}\label{ind2}
& E\exp\left(itA_\ell^{-1} (a_1X_{n_1}+\cdots +a_{k-1}X_{n_{k-1}}
+a_kX_{n} +a_{k+1}Y_{k+1} +\cdots + a_\ell Y_\ell)\right)\cr
&\qquad - E\exp \left( itA_\ell^{-1} (a_1X_{n_1}+\cdots + a_{k-1}X_{n_{k-1}} +a_kY_k +\cdots
+ a_\ell Y_\ell)\right)\\
& \le {{\varep}\over{2^k}} \psi
(a_1,\ldots ,a_\ell) \nonumber
\end{align}
for every $(a_1,\ldots,a_\ell)\in R^\ell$, $\ell >k$ and $n\ge n_k$.  This completes
the $k$-th induction step; the so constructed sequence $(n_k)$
obviously satisfies
\begin{align*}
& E\exp\left( itA_\ell^{-1} (a_1X_{n_1}+\cdots +a_\ell X_{n_\ell})\right)
- E\exp \left( itA_\ell ^{-1} (a_1Y_1 +\cdots + a_\ell Y_\ell)\right) \\
&\le \varep \psi (a_1,\ldots ,a_\ell)
\end{align*}
for every $\ell \ge 1$  and  $(a_1,\ldots, a_\ell)\in R^\ell$, i.e.\ (\ref{333})
is valid. Since in the $k$-th induction step $n_k$ was chosen in such a way that the corresponding
inequalities (\ref{ind1}) (for $k=1$) and (\ref{ind2}) (for $k>1$) hold not only for $n=n_k$, but
for all $n>n_k$ as well, relation (\ref{333}) remains valid for any further subsequence of $(X_{n_k})$.

To complete the proof of our theorem, it suffices to show that for any $t\in \mathbb{R}$
and any real sequence $(a_k)$ satisfying (\ref{an}) we have
\begin{equation}\label{chfu}
E\exp \left(it A_k^{-1} \sum_{j=1}^k a_j Y_j\right) \longrightarrow \exp(-c|t|^\alpha) \qquad \text{as} \ k\to\infty.
\end{equation}
Together with (\ref{333}) and the fact that (\ref{333}) remains valid for any further subsequence of $(X_{n_k})$
as well, this implies that for any $\ve>0$ and $t\in \mathbb{R}$ there exists an increasing sequence $(n_k)$ of
positive integers (depending on $\ve$ and $t$) such that for any further subsequence $(n_k')$ of $(n_k)$ we have
\begin{equation*}
\left| E\exp \left(it A_k^{-1} \sum_{j=1}^k a_j X_{n_j'}\right) - \exp(-c|t|^\alpha)\right| <\ve
\end{equation*}
for any $k\ge k_0 (\varepsilon, t)$ and any $(a_k)$ satisfying (\ref{an}). By a diagonal argument this shows that
there exists a sequence $(m_k)$ satisfying, together all of its subsequences, the relation
\begin{equation*}
E\exp \left(it A_k^{-1} \sum_{j=1}^k a_j X_{m_j}\right) \longrightarrow \exp(-c|t|^\alpha)
\end{equation*}
for any rational $t\in \mathbb{R}$ and any $(a_k)$ satisfying (\ref{an}), which implies that
$$
A_k^{-1} \sum_{j=1}^k a_j X_{m_j} \overset{d}{\longrightarrow} G_{\alpha, c},
$$
completing the ptoof of Theorem \ref{theorem2}. To verify (\ref{chfu}), let us note that conditionally
on $({\bf X}, \mu)$, $Y_j$ are i.i.d.\ with conditional characteristic function $\varphi$ satisfying (\ref{phi}),
which implies, in view of the remark after the proof of of Lemma \ref{lemma1}, that setting $S_k=\sum_{j=1}^k a_j Y_j$,
\begin{equation}\label {dom}
E\exp\left(itA_k^{-1}S_k | {\bf X}, \mu\right)\longrightarrow \exp(-c|t|^\alpha).
\end{equation}
Integrating the last relation and using the dominated convergence theorem we get (\ref{chfu}).

\end{document}